\documentclass[notitlepage,11pt]{article}
\catcode`\@=11
\@addtoreset{equation}{section}

\catcode`\@=12

\usepackage{latexsym}
\usepackage{amsmath}
\usepackage{amsfonts}
\usepackage{amssymb}
\usepackage{mathrsfs}
\usepackage{comment}

 \usepackage[colorlinks=true]{hyperref}
\hypersetup{urlcolor=blue, citecolor=red}

\newtheorem{theorem}{Theorem}[section]

\newtheorem{lemma}[theorem]{Lemma}

\newtheorem{remark}{Remark}

\def\R{\mathbb{R}}

\def\S{\mathbb{S}}

\def\rellich{\mu}
\def\hardy{\nu}

\def\proof{\noindent{\textbf{Proof. }}}
\def\QED{\hfill {$\square$}\goodbreak \medskip}

\begin{document}

\title{Radial and non radial ground states for a class of dilation invariant fourth order semilinear elliptic equations on $\R^{n}$}

\author{}



\date{}

{\Large
\centerline{Radial and non radial ground states for a class of dilation}
\centerline{invariant fourth order semilinear elliptic equations on $\R^{n}$}
}
\bigskip

\centerline{\scshape Paolo Caldiroli}
\medskip

{\footnotesize
 \centerline{Dipartimento di Matematica, Universit\`a di Torino}
 \centerline{via Carlo Alberto, 10 -- 10123 Torino, Italy}
 \centerline{email: \tt{paolo.caldiroli@unito.it}}
}
\bigskip

\begin{abstract}
We prove existence of extremal functions for some Rellich-Sobolev type inequalities involving the $L^{2}$ norm of the Laplacian as a leading term and the $L^{2}$ norm of the gradient, weighted with a Hardy potential. Moreover we exhibit a breaking symmetry phenomenon when the nonlinearity has a growth close to the critical one and the singular potential increases in strength.  
\vspace{6pt}\vspace{6pt}
\end{abstract}

\section{Introduction}
Recent years have seen a growing interest towards problems shaped on
\begin{equation}
\label{eq:pb0}
\left\{\begin{array}{l}\Delta^{2}u=|x|^{-\beta}|u|^{q-2}u\quad\text{on }\R^{n}\\ \int_{\R^{n}}|\Delta u|^{2}~\!dx<\infty
\end{array}\right.
\end{equation}
where $q>2$ and 
\begin{equation}
\label{eq:exponent}
\beta=n-\frac{q(n-4)}{2}~\!. 
\end{equation}
This choice of $\beta$ makes problems (\ref{eq:pb0}) invariant with respect to the action of the weighted dilation group
\begin{equation}
\label{eq:dilation}
\rho\mapsto u_{\rho}(x)=\rho^{\frac{n-4}{2}}u(\rho x)\quad(\rho>0).
\end{equation}
Problems (\ref{eq:pb0}) are variational in nature and rely on the validity of inequalities of the form
\begin{equation}
\label{eq:RS}
S_{q}\left(\int_{\R^{n}}|x|^{-\beta}|u|^{q}~\!dx\right)^{2/q}\le \int_{\R^{n}}|\Delta u|^{2}~\!dx\quad\forall u\in C^{\infty}_{c}(\R^{n}).
\end{equation}
In dimension $n\ge 5$ such inequalities can be obtained when 
$$
2\le q\le 2^{**}:=\frac{2n}{n-4}
$$
by interpolating the Rellich inequality (see \cite{Rel54}, \cite{Rel69})
$$
\left(\frac{n(n-4)}{4}\right)^{2}\int_{\R^{n}}|x|^{-4}|u|^{2}~\!dx\le \int_{\R^{n}}|\Delta u|^{2}~\!dx\quad\forall u\in C^{\infty}_{c}(\R^{n})
$$
with the Sobolev embedding
$$
S\left(\int_{\R^{n}}|u|^{2^{**}}~\!dx\right)^{2/2^{**}}\le \int_{\R^{n}}|\Delta u|^{2}~\!dx\quad\forall u\in C^{\infty}_{c}(\R^{n})~\!.
$$
We refer to the paper \cite{CaMu11} and to its bibliography for a deeper discussion on the inequalities (\ref{eq:RS}) and some generalizations. 

Problems (\ref{eq:pb0}) and many variants of them have been investigated in several works. Limiting ourselves to problems concerning entire solutions for equations ruled by the biharmonic operator or equivalent systems, we quote \cite{AlO}, \cite{BaMu}, \cite{Lions}, \cite{Mit93}, \cite{Mit96}, \cite{Mus13}, \cite{NousSwanYang} and the monography \cite{GhMo13} and we refer to the references therein contained. 

In this paper we study a variant of (\ref{eq:pb0}) characterized by the presence of a dilation invariant (hence, non compact) additional term containing lower order derivatives and whose shape preserves the variational character of the problem. More precisely we are interested in the existence of ground states for the problems
\begin{equation}
\label{eq:P}
\left\{
\begin{array}{ll}
\Delta^{2}u+\lambda~\!\mathrm{div}(|x|^{-2}\nabla u)=|x|^{-\beta}|u|^{q-2}u&\text{on }\R^{n}\\
\int_{\R^{n}}|\Delta u|^{2}~\!dx<\infty
\end{array}
\right.
\end{equation}
where $n\ge 5$, $q\in(2,2^{**}]$, $\lambda\in\R$ and $\beta$ is like in (\ref{eq:exponent}). The novelty with respect to the known literature stays in the term $\lambda~\!\mathrm{div}(|x|^{-2}\nabla u)$ containing the Hardy potential with strength ruled by the parameter $\lambda$.

By \emph{ground state} for (\ref{eq:P}) we mean a weak nontrivial solution of (\ref{eq:P}) belonging to the Sobolev space $D^{2,2}(\R^{n})$ and characterized as a minimum point for 
\begin{equation}
\label{eq:Sq}
S_{q}(\lambda):=\inf_{\scriptstyle u\in D^{2,2}(\R^{n})\atop\scriptstyle u\ne 0}\frac{\displaystyle\int_{\R^{n}}|\Delta u|^{2}~\!dx-\lambda\int_{\R^{n}}|x|^{-2}|\nabla u|^{2}~\!dx}{\displaystyle\left(\int_{\R^{n}}|x|^{-\beta}|u|^{q}~\!dx\right)^{{2}/{q}}}~\!.
\end{equation}
Here $D^{2,2}(\R^{n})$ is the space defined as the completion of $C^{\infty}_{c}(\R^{n})$ with respect to the norm
$$
\|u\|_{D^{2,2}}^{2}:=\sum_{i,j=1}^{n}\int_{\R^{n}}\left|\frac{\partial^{2}u}{\partial x_{i}\partial x_{j}}\right|^{2}~\!dx~\!.
$$
The restriction on the dimension $n\ge 5$ guarantees the Sobolev embedding for the space $D^{2,2}(\R^{n})$ into $L^{2^{**}}$.
In order to ensure that $S_{q}(\lambda)>0$ we have to take 
\begin{equation}
\label{eq:Lambda}
\lambda<\Lambda:=\inf_{\scriptstyle u\in D^{2,2}(\R^{n})\atop\scriptstyle u\ne 0}\frac{\displaystyle\int_{\R^{n}}|\Delta u|^{2}~\!dx}{\displaystyle\int_{\R^{n}}|x|^{-2}|\nabla u|^{2}~\!dx}~\!,
\end{equation}
$q\in[2,2^{**}]$, and $\beta$ as in (\ref{eq:exponent}). In \cite{TerZog} (see also \cite{AGS}, \cite{AS}, \cite{CaMu-Bangalore}, \cite{GhMo11} and \cite{Mor}) it was proved that if $n\ge 5$ then 
$$
\Lambda=\frac{n^{2}}{4}.
$$ 
We notice that, as well as (\ref{eq:pb0}), also problems (\ref{eq:P}) turn out to be invariant under the weighted dilation (\ref{eq:dilation}). As a consequence, the corresponding variational problems exhibit a lack of compactness. We can show the following result.

\begin{theorem}
\label{T:ground-state}
Let $n\ge 5$.
\begin{itemize}
\item[(i)]
For $q\in(2,2^{**})$ problem (\ref{eq:P}) admits a ground state for every $\lambda<\Lambda$.
\item[(ii)]
For $q=2^{**}$ problem (\ref{eq:P}) admits a ground state if and only if $0\le\lambda<\Lambda$. Moreover for every $\lambda\le 0$ the infimum $S_{2^{**}}(\lambda)$ equals the Sobolev constant of the embedding of $D^{2,2}(\R^{n})$ into $L^{2^{**}}$.
\end{itemize}
\end{theorem}

We can drop the upper bound on $q$ by looking for a \emph{radial ground state} for problem (\ref{eq:P}), namely, a non trivial, radial weak solution of (\ref{eq:P}) characterized as a minimum point for 
$$
S_{q}^{\mathrm{rad}}(\lambda):=\inf_{\scriptstyle u\in D^{2,2}_{\text{rad}}(\R^{n})\atop\scriptstyle u\ne 0}\frac{\displaystyle\int_{\R^{n}}|\Delta u|^{2}~\!dx-\lambda\int_{\R^{n}}|x|^{-2}|\nabla u|^{2}~\!dx}{\displaystyle\left(\int_{\R^{n}}|x|^{-\beta}|u|^{q}~\!dx\right)^{{2}/{q}}}
$$
where $D^{2,2}_{\text{rad}}(\R^{n})$ is the space of radial functions belonging to $D^{2,2}(\R^{n})$. We have that:

\begin{theorem}
\label{T:radial-ground-state}
If $n\ge 5$, $\lambda<\Lambda$ and $q\in(2,\infty)$ then problem (\ref{eq:P}) admits a radial ground state. Moreover such a ground state is positive and is unique up to the weighted dilation (\ref{eq:dilation}).
\end{theorem}

When $2<q\le 2^{**}$ we can compare the infima values $S_{q}(\lambda)$ and $S_{q}^{\mathrm{rad}}(\lambda)$ and in some cases we can observe a breaking symmetry phenomenon. More precisely we have a first result stated as follows.

\begin{theorem}
\label{T:breaking-symmetry-1}
For every $n\ge 5$ and $\lambda<0$ there exists $q_{\lambda,n}\in(2,2^{**})$ such that if $q\in(q_{\lambda,n},2^{**}]$ then $S_{q}(\lambda)<S_{q}^{\mathrm{rad}}(\lambda)$. In particular if $q\in(q_{\lambda,n},2^{**})$ then the ground state for problem (\ref{eq:P}) is non radial.
\end{theorem}

The previous result is obtained just by noticing that $S_{2^{**}}(\lambda)<S_{2^{**}}^{\mathrm{rad}}(\lambda)$ and using the continuity of the mappings $q\mapsto S_{q}(\lambda)$ and $q\mapsto S_{q}^{\mathrm{rad}}(\lambda)$. We have no information on $q_{\lambda,n}$, that is, on the range of $q$'s for which breaking symmetry occurs. More precise estimates are stated in the next theorem.

\begin{theorem}
\label{T:breaking-symmetry-2}
If $q\in(2,2^{**})$ satisfies
\begin{equation}
\label{eq:BS}
\frac{3q+2}{q(q-2)}<\frac{(n-4)^{2}}{n-1}
\end{equation}
then for $\lambda<0$ with $|\lambda|$ large enough (depending on $q$) one has that $S_{q}(\lambda)<S_{q}^{\mathrm{rad}}(\lambda)$. In such a case problem (\ref{eq:P}) admits at least two non trivial solutions and its ground state is non radial.
\end{theorem}

Condition (\ref{eq:BS}) is fulfilled if $q$ is close to the critical exponent $2^{**}$. More precisely, setting
$$
q_{n}=1+a_{n}+\sqrt{(1+a_{n})^{2}+\tfrac{4}{3}a_{n}}\quad\text{where}\quad a_{n}=\frac{3(n-1)}{2(n-4)^{2}}~\!,
$$ 
one has that $q_{n}>2$ and (\ref{eq:BS}) holds true for $q\in(q_{n},2^{**})$. However the interval $(q_{n},2^{**})$ is nonempty just for $n\ge 7$. On the other hand we observe that $q_{n}\to 2$ (as well as $2^{**}$) as $n\to\infty$. We also notice that in the linear case, namely when $q=2$, we have that $S_{2}(\lambda)=S_{2}^{\mathrm{rad}}(\lambda)$ for every $\lambda<\Lambda$ (see Remark \ref{R:linear}).

We point out that similar, actually sharper, existence results of radial and non radial ground states, as well as breaking symmetry, have been proved in \cite{BaMu} for the problem
\begin{equation}
\label{eq:P1}
\left\{
\begin{array}{ll}
\Delta^{2}u-\lambda|x|^{-4}u=|x|^{-\beta}|u|^{q-2}u&\text{on }\R^{n}\\
\int_{\R^{n}}|\Delta u|^{2}~\!dx<\infty~\!.
\end{array}
\right.
\end{equation}
In fact, in the critical case $q=2^{**}$ (and $\beta=0$) problems (\ref{eq:P}) and (\ref{eq:P1}) can be viewed as higher order versions of the problem
\begin{equation}
\label{eq:P2}
\left\{
\begin{array}{ll}
-\Delta u-\lambda|x|^{-2}u=|u|^{2^{*}-2}u&\text{on }\R^{n}\\
\int_{\R^{n}}|\nabla u|^{2}~\!dx<\infty
\end{array}
\right.
\end{equation}
where $2^{*}=2n/(n-2)$ and $n\ge 3$. As proved in \cite{Ter96}, also (\ref{eq:P2}) admits a couple of non trivial solutions, characterized as radial and non radial ground states, and they are different when $-\lambda<\lambda_{0}$ for some $\lambda_{0}<0$. See also \cite{CatrinaWang} for a class of second order problems generalizing (\ref{eq:P2}) and displaying breaking symmetry.  

Indeed problems (\ref{eq:P}), (\ref{eq:P1}) and (\ref{eq:P2}) share similar features: all of them are based on a suitable functional inequality and their solutions can be found as extremal functions for such inequality. Moreover, roughly speaking, the breaking symmetry is due to the fact that, as $\lambda\to-\infty$, the singular potential in the corresponding lower order term becomes more and more important and changes the topology of lower sublevel sets.

\section{Proof of Theorem \ref{T:ground-state}}

A key tool in our argument is the following compactness lemma. This result is an adaptation of a tool already used in previous works, like \cite{BaMu} or \cite{CaMu11}.

\begin{lemma}
\label{L:compact}
Let $R>0$ and let $(u_{k})$ be a sequence in $D^{2,2}(\R^{n})$ satisfying
\begin{gather}
\label{eq:weak-conv}
u_{k}\to 0\text{ weakly in }D^{2,2}(\R^{n})\\
\label{eq:almost-solution}
\Delta^{2}u_{k}+\lambda~\!\mathrm{div}(|x|^{-2}\nabla u_{k})-|x|^{-\beta}|u_{k}|^{q-2}u_{k}\to 0\text{ in $(D^{2,2}(\R^{n}))'$}\\
\label{eq:small}
\limsup\int_{B_{R}}|x|^{-\beta}|u_{k}|^{q}~\!dx<S_{q}(\lambda)^{q/(q-2)}.
\end{gather}
Then $|x|^{-\beta}|u_{k}|^{q}\to 0$ strongly in $L^{1}_{\mathrm{loc}}(B_{R})$.
\end{lemma}

\proof
Fix $R'\in(0,R)$ and take a cut-off function $\varphi\in C^{\infty}_{c}(B_{R})$ such that $\varphi=1$ on $B_{R'}$. We point out that the sequence $(\varphi^{2} u_{k})$ is bounded in $D^{2,2}(\R^{n})$. Using $\varphi^{2} u_{k}$ as a test function in (\ref{eq:almost-solution}) we obtain
\begin{equation}
\label{eq:almost-sol-test}
\int_{\R^{n}}\varphi^{2} u_{k}\Delta^{2}u_{k}~\!dx-\lambda\int_{\R^{n}}|x|^{-2}\nabla(\varphi^{2} u_{k})\cdot\nabla u_{k}~\!dx=\int_{\R^{n}}|x|^{-\beta}\varphi^{2}|u_{k}|^{q}~\!dx+o(1).
\end{equation}
By (\ref{eq:weak-conv}) $u_{k}\to 0$ weakly in $H^{2}_{\mathrm{loc}}(\R^{n})$ and then, by compactness, $u_{k}\to 0$ strongly in $H^{1}(B_{R})$. Hence we have that
\begin{gather*}
\int_{\R^{n}}|\Delta(\varphi u_{k})|^{2}~\!dx=\int_{\R^{n}}\varphi^{2}|\Delta u_{k}|^{2}~\!dx+o(1)\\
\int_{\R^{n}}(\Delta u_{k})\Delta(\varphi^{2} u_{k})~\!dx=\int_{\R^{n}}\varphi^{2}|\Delta u_{k}|^{2}~\!dx+o(1)\\
\int_{\R^{n}}|x|^{-2}\nabla u_{k}\cdot\nabla(\varphi^{2}u_{k})~\!dx=
\int_{\R^{n}}|x|^{-2}|\nabla(\varphi u_{k})|^{2}~\!dx+o(1)
\end{gather*}
Then, after integration by parts,
\begin{gather*}
\label{eq:Delta2}
\int_{\R^{n}}\varphi^{2}u_{k}\Delta^{2}u_{k}~\!dx=\int_{\R^{n}}|\Delta(\varphi u_{k})|^{2}~\!dx+o(1)\\
\label{eq:nabla}
\int_{\R^{n}}\varphi^{2}u_{k}~\!\mathrm{div}(|x|^{-2}\nabla u_{k})~\!dx=-\int_{\R^{n}}|x|^{-2}|\nabla(\varphi u_{k})|^{2}~\!dx+o(1).
\end{gather*}
Consequently (\ref{eq:almost-sol-test}) reduces to
\begin{equation}
\label{eq:almost-sol-test2}
\int_{\R^{n}}|\Delta(\varphi u_{k})|^{2}~\!dx-\lambda\int_{\R^{n}}|x|^{-2}|\nabla(\varphi u_{k})|^{2}~\!dx=\int_{\R^{n}}|x|^{-\beta}\varphi^{2}|u_{k}|^{q}~\!dx+o(1).
\end{equation}
By (\ref{eq:small}) there exists $\varepsilon_{0}>0$ such that
\begin{equation}
\label{eq:epsilon}
\int_{\R^{n}}|x|^{-\beta}|\varphi u_{k}|^{q}~\!dx\le\varepsilon_{0}<S_{q}(\lambda)^{q/(q-2)}\quad\forall k\text{ large.}
\end{equation}
Therefore, using the H\"{o}lder inequality and (\ref{eq:epsilon}), we estimate
\begin{equation}
\label{eq:almost-sol-test3}
\int_{\R^{n}}|x|^{-\beta}\varphi^{2}|u_{k}|^{q}~\!dx\le \varepsilon_{0}^{(q-2)/q}\left(\int_{\R^{n}}|x|^{-\beta}|\varphi u_{k}|^{q}~\!dx\right)^{2/q}.
\end{equation}
On the other side, by definition of $S_{q}(\lambda)$,
\begin{equation}
\label{eq:almost-sol-test4}
\int_{\R^{n}}|\Delta(\varphi u_{k})|^{2}~\!dx-\lambda\int_{\R^{n}}|x|^{-2}|\nabla(\varphi u_{k})|^{2}~\!dx\ge S_{q}(\lambda)\left(\int_{\R^{n}}|x|^{-\beta}|\varphi u_{k}|^{q}~\!dx\right)^{2/q}.
\end{equation}
Therefore from (\ref{eq:almost-sol-test2})--(\ref{eq:almost-sol-test4}) it follows that
$$
S_{q}(\lambda)\left(\int_{\R^{n}}|x|^{-\beta}|\varphi u_{k}|^{q}~\!dx\right)^{2/q}\le\varepsilon_{0}^{(q-2)/q}\left(\int_{\R^{n}}|x|^{-\beta}|\varphi u_{k}|^{q}~\!dx\right)^{2/q}+o(1).
$$
As $\varepsilon_{0}<S_{q}(\lambda)^{q/(q-2)}$ we infer that 
$$
\int_{\R^{n}}|x|^{-\beta}|\varphi u_{k}|^{q}~\!dx\to 0
$$
and then, since $\varphi=1$ on $B_{R'}$ and $R'$ is arbitrary in $(0,R)$, $|x|^{-\beta}|u_{k}|^{q}\to 0$ strongly in $L^{1}_{\mathrm{loc}}(B_{R})$.
\QED

Now let us proceed with the proof of Theorem \ref{T:ground-state}.
Using Ekeland's variational principle (see \cite{Str} Chapt.~1, Sect.~5) and the invariance under the weighted dilation (\ref{eq:dilation}) we can find a minimizing sequence $(u_{k})\subset D^{2,2}(\R^{n})$ for problem (\ref{eq:Sq}), satisfying (\ref{eq:almost-solution}) and
\begin{gather}
\nonumber
\int_{\R^{n}}\left(|\Delta u_{k}|^2-\lambda|x|^{-2}|\nabla u_{k}|^{2}\right)~\!dx=S_{q}(\lambda)^{{q}/{q-2}}+o(1)\\
\label{eq:q-norm}
\int_{\R^{n}}|x|^{-\beta}|u_{k}|^{q}~\!dx=S_{q}(\lambda)^{{q}/{q-2}}+o(1)\\
\label{eq:B2}
\int_{B_{2}}|x|^{-\beta}|u_{k}|^{q}~\!dx=\frac{1}{2}S_{q}(\lambda)^{{q}/{q-2}}.
\end{gather}
Since $\lambda<\Lambda$, we have that
$$
\sup\|\Delta u_{k}\|_{L^{2}(\R^{n})}<\infty.
$$
It is known that for $n\ge 5$ in the space $D^{2,2}(\R^{n})$ the $L^{2}$-norm of the laplacian is equivalent to the $D^{2,2}$-norm (see \cite{Lin86} and Remark 2.3 in \cite{CaMu11}). Hence the sequence $(u_{k})$ is bounded in $D^{2,2}(\R^{n})$ and then it admits a subsequence, still denoted $(u_{k})$, weakly converging to some $u\in D^{2,2}(\R^{n})$. If $u\ne 0$, then $u$ is a minimizer for $S_{q}(\lambda)$ and $u_{k}\to u$ strongly in $D^{2,2}(\R^{n})$. The proof of this fact is definitely standard: one can adapt to our situation a well known argument (see, e.g., \cite{Str}, Chapt.~1, Sect.~4). Hence we have to exclude that $u=0$. We argue by contradiction, assuming that $u=0$. In this case, by Lemma \ref{L:compact}
$$
\int_{B_{1}}|x|^{-\beta}|u_{k}|^{q}\!~dx=o(1).
$$
Therefore, by (\ref{eq:B2}),
\begin{equation}
\label{eq:annulus}
\int_{B_{2}\setminus B_{1}}|x|^{-\beta}|u_{k}|^{q}\!~dx\to\frac{1}{2}S_{q}(\lambda)^{{q}/{q-2}}.
\end{equation}
Let us distinguish the cases of subcritical or critical exponent.\medskip

\noindent
$(i)$ If $q\in(2,2^{**})$, since $u_{k}\to 0$ weakly in $H^{2}_{\mathrm{loc}}(\R^{n})$, the Rellich compactness Theorem implies that $u_{k}\to 0$ strongly in $L^{q}(B_{2}\setminus B_{1})$, contradicting (\ref{eq:annulus}). Hence in this case the weak limit $u$ cannot be zero and the proof is complete.
\medskip

\noindent
$(ii)$ Now we study the case of critical exponent $q=2^{**}$. Notice that for such a value of $q$, the corresponding exponent $\beta$ is null. Firstly we take $\lambda\in(0,\Lambda)$. Let us fix a cut-off function $\varphi\in C^{\infty}_{c}(\R^{n}\setminus\{0\})$ such that $0\le\varphi\le 1$ and $\varphi(x)=1$ for $1\le|x|\le 2$. Arguing as in the first part of the proof of Lemma \ref{L:compact} we obtain that 
\begin{equation}
\label{eq:almost-sol-crit1}
\int_{\R^{n}}|\Delta(\varphi u_{k})|^{2}~\!dx-\lambda\int_{\R^{n}}|x|^{-2}|\nabla(\varphi u_{k})|^{2}~\!dx=\int_{\R^{n}}\varphi^{2}|u_{k}|^{2^{**}}~\!dx+o(1).
\end{equation}
Since
$$
\int_{\R^{n}}|x|^{-2}|\nabla(\varphi u_{k})|^{2}~\!dx\le C\int_{B_{2}\setminus B_{1}}\left(|\nabla u_{k}|^{2}+u_{k}^{2}\right)~\!dx
$$
and $u_{k}\to 0$ strongly in $H^{1}_{\mathrm{loc}}(\R^{n})$, (\ref{eq:almost-sol-crit1}) reduces to
\begin{equation}
\label{eq:almost-sol-crit2}
\int_{\R^{n}}|\Delta(\varphi u_{k})|^{2}~\!dx=\int_{\R^{n}}\varphi^{2}|u_{k}|^{2^{**}}~\!dx+o(1).
\end{equation}
On one side, using the H\"{o}lder inequality and (\ref{eq:q-norm}), we estimate
$$
\int_{\R^{n}}\varphi^{2}|u_{k}|^{2^{**}}~\!dx\le S_{2^{**}}(\lambda)\left(\int_{\R^{n}}|\varphi u_{k}|^{2^{**}}~\!dx\right)^{2/2^{**}}+o(1).
$$
On the other side
$$
\int_{\R^{n}}|\Delta(\varphi u_{k})|^{2}~\!dx\ge S_{2^{**}}(0)\left(\int_{\R^{n}}|\varphi u_{k}|^{2^{**}}~\!dx\right)^{2/2^{**}}
$$
Therefore from (\ref{eq:almost-sol-crit2}) it follows that
\begin{equation}
\label{eq:almost-sol-crit3}
S_{2^{**}}(0)\left(\int_{\R^{n}}|\varphi u_{k}|^{2^{**}}~\!dx\right)^{2/2^{**}}\le S_{2^{**}}(\lambda)\left(\int_{\R^{n}}|\varphi u_{k}|^{2^{**}}~\!dx\right)^{2/2^{**}}+o(1).
\end{equation}
Now we claim that
\begin{equation}
\label{eq:strict-inequality}
S_{2^{**}}(\lambda)<S_{2^{**}}(0)\quad\forall\lambda\in(0,\Lambda).
\end{equation}
Indeed, let us notice that $S_{2^{**}}(0)$ equals the Sobolev constant
$$
S^{**}:=\inf_{\scriptstyle u\in D^{2,2}(\R^{n})\atop\scriptstyle u\ne 0}\frac{\displaystyle\int_{\R^{n}}|\Delta u|^{2}~\!dx}{\displaystyle\left(\int_{\R^{n}}|u|^{2^{**}}dx\right)^{2/2^{**}}}~.
$$
It is known that $S^{**}$ is achieved in $D^{2,2}(\R^{n})$ by $U(x)=(1+|x|^{2})^{(4-n)/2}$ (see for instance \cite{Sw}). Hence, as $\lambda>0$,
$$
S_{2^{**}}(0)=\frac{\displaystyle\int_{\R^{n}}|\Delta U|^{2}~\!dx}{\displaystyle\left(\int_{\R^{n}}|U|^{2^{**}}dx\right)^{2/2^{**}}}>
\frac{\displaystyle\int_{\R^{n}}|\Delta U|^{2}~\!dx-\lambda\int_{\R^{n}}|x|^{-2}|\nabla U|^{2}~\!dx}{\displaystyle\left(\int_{\R^{n}}|U|^{2^{**}}dx\right)^{2/2^{**}}}\ge S_{2^{**}}(\lambda)
$$
and then (\ref{eq:strict-inequality}) holds. Finally (\ref{eq:almost-sol-crit3}) and (\ref{eq:strict-inequality}) imply that $u_{k}\to 0$ strongly in $L^{2^{**}}(B_{2}\setminus B_{1})$, in contradiction with (\ref{eq:annulus}). Hence we proved that if $\lambda\in(0,\Lambda)$ then $S_{2^{**}}(\lambda)$ is achieved. In a standard way one shows the existence of a ground state. The case $\lambda=0$ is known, as mentioned before. It remains to study the case $\lambda<0$. Firstly we show that if $\lambda<0$ then 
\begin{equation}
\label{eq:S-equality}
S_{2^{**}}(\lambda)=S_{2^{**}}(0)\!~.
\end{equation}
Indeed it is clear that $S_{2^{**}}(\lambda)\ge S_{2^{**}}(0)$. Let us check the opposite inequality: for every $u\in C^{\infty}_{c}(\R^{n})\setminus\{0\}$ we set $u_{y}(x)=u(x-y)$ and we estimate
\begin{equation*}
\begin{split}
S_{2^{**}}(\lambda)&\le\lim_{|y|\to\infty}\frac{\displaystyle\int_{\R^{n}}|\Delta u_{y}|^{2}~\!dx-\lambda\int_{\R^{n}}|x|^{-2}|\nabla u_{y}|^{2}~\!dx}{\displaystyle\left(\int_{\R^{n}}|u_{y}|^{2^{**}}dx\right)^{2/2^{**}}}\\
&=\lim_{|y|\to\infty}\frac{\displaystyle\int_{\R^{n}}|\Delta u|^{2}~\!dx-\lambda\int_{\R^{n}}|x+y|^{-2}|\nabla u|^{2}~\!dx}{\displaystyle\left(\int_{\R^{n}}|u|^{2^{**}}dx\right)^{2/2^{**}}}=\frac{\displaystyle\int_{\R^{n}}|\Delta u|^{2}~\!dx}{\displaystyle\left(\int_{\R^{n}}|u|^{2^{**}}dx\right)^{2/2^{**}}}\!~.
\end{split}
\end{equation*}
By the arbitrariness of $u\in C^{\infty}_{c}(\R^{n})\setminus\{0\}$ and since $C^{\infty}_{c}(\R^{n})$ is dense in $D^{2,2}(\R^{n})$ we obtain that $S_{2^{**}}(\lambda)\le S_{2^{**}}(0)$. Hence (\ref{eq:S-equality}) holds. Moreover if $\lambda<0$ the infimum $S_{2^{**}}(\lambda)$ cannot be achieved. Otherwise, if $u\in D^{2,2}(\R^{n})$ would be a minimizer for $S_{2^{**}}(\lambda)$, then it would be a minimizer also for $S_{2^{**}}(0)$. In particular $\int_{\R^{n}}|x|^{-2}|\nabla u|^{2}~\!dx=0$, that is $u=0$, which is impossible. Thus the proof is complete.
\QED

\section{Proof of Theorem \ref{T:radial-ground-state}}
Let us introduce the Emden-Fowler transform, defined as follows: for every $u\in D^{2,2}_{\mathrm{rad}}(\R^{n})$ let $w\colon(0,\infty)\to\R$ be such that
\begin{equation}
\label{eq:EF}
u(x)=|x|^{\frac{4-n}{2}}w(-\log|x|).
\end{equation}
\begin{lemma}
\label{L:EF}
For $n\ge 5$ the mapping $u\mapsto w=Tu$ defines an isomorphism between $D^{2,2}_{\mathrm{rad}}(\R^{n})$ and $H^{2}(\R)$. Moreover, setting $\omega_{n}=|\mathbb{S}^{n-1}|$, one has that
\begin{equation*}
\begin{split}
&\int_{\R^{n}}|\Delta u|^{2}~\!dx=\omega_{n}\int_{\R}\left(|w''|^{2}+
2\left(\rellich_{n}+2\right)|w'|^{2}+\rellich_{n}^{2}|w|^{2}\right)~\!dt\\
&\int_{\R^{n}}|x|^{-2}|\nabla u|^{2}~\!dx=\omega_{n}\int_{\R}\left(|w'|^{2}+\hardy_{n}|w|^{2}\right)~\!dt\\
&\int_{\R^{n}}|x|^{-\beta}|u|^{q}~\!dx=\omega_{n}\int_{\R}|w|^{q}~\!dt
\end{split}
\end{equation*}
where 
\begin{equation}
\label{eq:rellich-hardy}
\rellich_{n}=\frac{n(n-4)}{4}\quad\text{and}\quad\hardy_{n}=\frac{(n-4)^{2}}{4}~\!.
\end{equation}
\end{lemma}
For the proof we refer to \cite{CaMu11}. In view of Lemma \ref{L:EF} we have that
\begin{equation}
\label{eq:radial-minimization}
S_{q}^{\mathrm{rad}}(\lambda)=\omega_{n}^{\frac{q-2}{q}}\inf_{\scriptstyle w\in H^{2}(\R)\atop\scriptstyle w\ne 0}\frac{\displaystyle\int_{\R}\left(|w''|^{2}+2a_{\lambda}|w'|^{2}+b_{\lambda}|w|^{2}\right)~\!dt}{\displaystyle\left(\int_{\R}|w|^{q}~\!dt\right)^{{2}/{q}}}
\end{equation}
where
\begin{equation}
\label{eq:alambda-blambda}
2a_{\lambda}=2\left(\rellich_{n}+2\right)-\lambda\quad\text{and}\quad b_{\lambda}=(\Lambda-\lambda)\hardy_{n}~\!.
\end{equation}
We point out that, thanks to the assumption $\lambda<\Lambda$, the values $a_{\lambda}$ and $b_{\lambda}$ are positive. Now we use the following key result, proved in \cite{BaMu}:

\begin{theorem}
\label{T:BM}
For every $a,b>0$ and $q>2$ the minimization problem
$$
\inf_{\scriptstyle w\in H^{2}(\R)\atop\scriptstyle w\ne 0}\frac{\displaystyle\int_{\R}\left(|w''|^{2}+2a|w'|^{2}+b|w|^{2}\right)~\!dt}{\displaystyle\left(\int_{\R}|w|^{q}~\!dt\right)^{{2}/{q}}}
$$
admits a minimum point. In addition, if $a^{2}\ge b$ then the minimum point is positive and unique, up to the natural invariances of the problem (i.e., translation, inversion, multiplication by a non zero constant).
\end{theorem}

In the case in consideration 
$$
a_{\lambda}^{2}-b_{\lambda}=\left(\frac{\lambda}{2}-(n-2)\right)^{2}.
$$
Hence by Theorem \ref{T:BM} there exists a positive function $w\in H^{2}(\R)$ which is a minimizer for the problem defined by the right hand side of (\ref{eq:radial-minimization}). Such a minimizer is unique up to translation, inversion, and multiplication by a non zero constant. Then, using Lemma \ref{L:EF}, we infer that the mapping $u$ defined by (\ref{eq:EF}) belongs to $D^{2,2}_{\mathrm{rad}}(\R^{n})$, is a positive minimizer for $S_{q}^{\mathrm{rad}}(\lambda)$ and is the unique minimizer up to the weighted dilation (\ref{eq:dilation}) and to a multiplicative constant. In a standard way one also infers that for a suitable $\alpha>0$ the mapping $\alpha u$ is a radial ground state for problem (\ref{eq:P}).
\QED

\section{Proof of Theorem \ref{T:breaking-symmetry-1}}

If $q=2^{**}$ and $\lambda<0$ then by Theorems \ref{T:ground-state} and \ref{T:radial-ground-state} $S_{2^{**}}(\lambda)=S^{**}$ is not attained whereas $S_{2^{**}}^{\mathrm{rad}}(\lambda)$ is so. Hence $S_{2^{**}}^{\mathrm{rad}}(\lambda)>S_{2^{**}}(\lambda)$. The mapping $q\mapsto S_{q}^{\mathrm{rad}}(\lambda)$ and $q\mapsto S_{q}(\lambda)$ are continuous (see Remark B.4 and Lemma B.5 in \cite{CaMu11}). Thus the conclusion readily follows.
\QED

\section{Proof of Theorem \ref{T:breaking-symmetry-2}}

Let us fix $q\in(2,2^{**})$. By Theorems \ref{T:ground-state} and \ref{T:radial-ground-state} for every $\lambda<\Lambda$ the infima $S_{q}(\lambda)$ and $S_{q}^{\mathrm{rad}}(\lambda)$ are attained in $D^{2,2}(\R^{n})$. Clearly $S_{q}(\lambda)\le S_{q}^{\mathrm{rad}}(\lambda)$. Assume that equality holds. Let $u$ be a radial minimizer of the functional
$$
J(v)=\frac{\displaystyle\int_{\R^{n}}|\Delta v|^{2}~\!dx-\lambda\int_{\R^{n}}|x|^{-2}|\nabla v|^{2}~\!dx}{\displaystyle\left(\int_{\R^{n}}|x|^{-\beta}|v|^{q}~\!dx\right)^{{2}/{q}}}\quad(v\in D^{2,2}(\R^{n})\setminus\{0\}).
$$
Since $J$ is homogeneous, i.e., $J(\alpha v)=J(v)$ for any $\alpha\ne 0$, we can assume that 
\begin{equation}
\label{eq:u-Lq}
\int_{\R^{n}}|x|^{-\beta}|u|^{q}~\!dx=1.
\end{equation}
Since $u$ is a minimizer for $J$ on the whole space $D^{2,2}(\R^{n})\setminus\{0\}$ we have that $J'(u)[v]=0$ and $J''(u)[v,v]\ge 0$ for every $v\in D^{2,2}(\R^{n})$. Setting
$$
A(v)=\int_{\R^{n}}|\Delta v|^{2}~\!dx-\lambda\int_{\R^{n}}|x|^{-2}|\nabla v|^{2}~\!dx\text{~~and~~}B(v)=\left(\int_{\R^{n}}|x|^{-\beta}|v|^{q}~\!dx\right)^{{2}/{q}},
$$
since $J'(u)=0$, $A''(u)[v,v]=2A(v)$ and $B(u)=1$, the condition $J''(u)[v,v]\ge 0$ reads
\begin{equation}
\label{eq:second-derivative}
A(v)\ge \frac{1}{2}A(u)B''(u)[v,v]\quad\text{for every }v\in D^{2,2}(\R^{n}).
\end{equation}
We take $v=u\varphi$ where $\varphi\in C^{\infty}(\S^{n-1})$ is an eigenfunction of the Laplace-Beltrami operator on the sphere $\S^{n-1}$ corresponding to the first positive eigenvalue $n-1$. We normalize $\varphi$ according to the condition
\begin{equation}
\label{eq:phi-L2}
\int_{\S^{n-1}}\varphi^{2}~\!d\sigma=|\S^{n-1}|.
\end{equation}
Moreover we notice that
\begin{equation}
\label{eq:phi-average}
\int_{\S^{n-1}}\varphi~\!d\sigma=0.
\end{equation}
Using (\ref{eq:u-Lq}), (\ref{eq:phi-L2}) and (\ref{eq:phi-average}), one computes $B''(u)[u\varphi,u\varphi]=2(q-1)$. Hence from (\ref{eq:second-derivative}) it follows that
\begin{equation}
\label{eq:second-derivative-bis}
A(u\varphi)\ge(q-1)A(u).
\end{equation}
In order to evaluate $A(u\varphi)$ we observe that 
\begin{equation*}
\begin{split}
&|\nabla(u\varphi)|^{2}=|\nabla u|^{2}\varphi^{2}+|x|^{-2}u^{2}|\nabla_{\!\sigma}\varphi|^{2}\\
&\Delta(u\varphi)=(\Delta u)\varphi+|x|^{-2}u\Delta_{\sigma}\varphi=(\Delta u-(n-1)|x|^{-2}u)\varphi.
\end{split}
\end{equation*}
Thanks to the above expressions, and since 
$$
\int_{\S^{n-1}}|\nabla_{\!\sigma}\varphi|^{2}~\!d\sigma=(n-1)\int_{\S^{n-1}}\varphi^{2}~\!d\sigma~\!,
$$
we infer that
\begin{equation}
\label{eq:A}
A(u\varphi)=A(u)+(n-1)(n-1-\lambda)\int_{\R^{n}}|x|^{-4}u^{2}~\!dx-2(n-1)\int_{\R^{n}}|x|^{-2}u\Delta u~\!dx~\!.
\end{equation}
Let us introduce the following shortened notation:
$$
U_{0}=\int_{\R^{n}}|x|^{-4}u^{2}~\!dx~\!,\quad U_{1}=\int_{\R^{n}}|x|^{-2}|\nabla u|^{2}~\!dx~\!,\quad U_{2}=\int_{\R^{n}}|\Delta u|^{2}~\!dx~\!.
$$
Since 
$$
\int_{\R^{n}}|x|^{-2}u\Delta u~\!dx=-(n-4)U_{0}-U_{1}~\!,
$$
from (\ref{eq:A}) and (\ref{eq:second-derivative-bis}) we obtain
\begin{equation}
\label{eq:second-derivative-ter}
U_{2}\le\left(\lambda+\frac{2(n-1)}{q-2}\right)U_{1}+\frac{(n-1)(3n-9-\lambda)}{q-2}~\!U_{0}~\!.
\end{equation}
Since $u$ is a radial minimizer for $S_{q}(\lambda)$, and is normalized according to the condition (\ref{eq:u-Lq}), the map $w$ defined by (\ref{eq:EF}) solves:
$$
w''''-2a_{\lambda}w''+b_{\lambda}w=S_{q}(\lambda)w^{q-1}
$$
where $a_{\lambda}$ and $b_{\lambda}$ are defined as in (\ref{eq:alambda-blambda}). By the energy conservation and since $w\in H^{2}(\R)$ one has that
$$
-w'''w'+\frac{1}{2}|w''|^{2}+a_{\lambda}|w'|^{2}-\frac{b_{\lambda}}{2}|w|^{2}+\frac{S_{q}(\lambda)}{q}|w|^{q}=0
$$
and then, after an integration by parts,
$$
3\int_{\R}|w''|^{2}~\!dt+2a_{\lambda}\int_{\R}|w'|^{2}~\!dt-b_{\lambda}\int_{\R}|w|^{2}~\!dt+\frac{2}{q}~\!S_{q}(\lambda)\int_{\R}|w|^{q}~\!dt=0~\!.
$$
By Lemma \ref{eq:EF}, by (\ref{eq:u-Lq}) and since $A(u)=S_{q}(\lambda)$, the previous equation can be written in the form:
\begin{equation}
\label{eq:u-L2}
\left(3+\frac{2}{q}\right)U_{2}
=\left(4(\rellich_{n}+2)+\lambda~\!\frac{q-2}{q}\right)U_{1}-2\hardy_{n}\left(2(\rellich_{n}+2)-\hardy_{n}+\lambda\right)U_{0}
\end{equation}
with $\rellich_{n}$ and $\hardy_{n}$ defined in (\ref{eq:rellich-hardy}). Then (\ref{eq:second-derivative-ter}) and (\ref{eq:u-L2}) imply 
\begin{equation*}
\begin{split}
&\left[-\frac{2q}{3q+2}~\!\lambda+\left(\frac{4q(\rellich_{n}+2)}{3q+2}-\frac{2(n-1)}{q-2}\right)\right]U_{1}\\
&\qquad\qquad\le\left[\left(\frac{2q\hardy_{n}}{3q+2}-\frac{n-1}{q-2}\right)\lambda+\left(-\frac{4q\hardy_{n}(n-2)}{3q+2}+\frac{3(n-1)(n-3)}{q-2}\right)\right]U_{0}.
\end{split}
\end{equation*}
Now we apply the Hardy inequality
\begin{equation}
\label{hardy-inequality}
\int_{\R^{n}}|x|^{-2}|\nabla v|^{2}~\!dx>\nu_{n}\int_{\R^{n}}|x|^{-4}v^{2}~\!dx
\end{equation}
to the mapping $u$, getting that $U_{1}>\hardy_{n}U_{0}$. Hence, taking $\lambda<0$ such that
\begin{equation}
\label{eq:lambda-condition}
-\frac{2q}{3q+2}~\!\lambda+\left(\frac{4q(\rellich_{n}+2)}{3q+2}-\frac{2(n-1)}{q-2}\right)>0
\end{equation}
we obtain
\begin{equation*}
\begin{split}
&-\frac{2q\hardy_{n}}{3q+2}~\!\lambda+\left(\frac{4q(\rellich_{n}+2)}{3q+2}-\frac{2(n-1)}{q-2}\right)\hardy_{n}\\
&\qquad\qquad\le\left(\frac{2q\hardy_{n}}{3q+2}-\frac{n-1}{q-2}\right)\lambda+\left(-\frac{4q\hardy_{n}(n-2)}{3q+2}+\frac{3(n-1)(n-3)}{q-2}\right).
\end{split}
\end{equation*}
Divinding by $-\lambda$ and letting $\lambda\to-\infty$ we get
\begin{equation}
\label{eq:q-condition}
\frac{2q\hardy_{n}}{3q+2}\le-\frac{2q\hardy_{n}}{3q+2}+\frac{n-1}{q-2}.
\end{equation}
Hence we proved that for $\lambda<0$ and $q\in(2,2^{**})$ a necessary condition in order that $S_{q}(\lambda)=S_{q}^{\mathrm{rad}}(\lambda)$ is that (\ref{eq:lambda-condition}) and (\ref{eq:q-condition}) hold. Since (\ref{eq:BS}) is equivalent to the failure of (\ref{eq:q-condition}), the thesis is proved.
\QED

\begin{remark}
\label{R:linear}
As far as concerns the minimization problems defined by $S_{2}(\lambda)$ and $S_{2}^{\mathrm{rad}}(\lambda)$, we can easily show that
\begin{equation}
\label{eq:S2}
S_{2}(\lambda)=S_{2}^{\mathrm{rad}}(\lambda)=\rellich_{n}^{2}-\lambda\hardy_{n}\quad\forall\lambda<\Lambda~\!.
\end{equation}
Indeed, using (\ref{eq:Lambda}) and the Hardy inequality (\ref{hardy-inequality}), for $\lambda<\Lambda$ we readily have that $S_{2}(\lambda)\ge(\Lambda-\lambda)\hardy_{n}$. Moreover, by Lemma \ref{L:EF},  we have that
$$
S_{2}^{\mathrm{rad}}(\lambda)=\inf_{\scriptstyle w\in H^{2}(\R)\atop\scriptstyle w\ne 0}\frac{\displaystyle\int_{\R}\left(|w''|^{2}+2a_{\lambda}|w'|^{2}+b_{\lambda}|w|^{2}\right)~\!dt}{\displaystyle\int_{\R}|w|^{2}~\!dt}
$$
with $a_{\lambda}$ and $b_{\lambda}$ as in (\ref{eq:alambda-blambda}). A standard dilation argument shows that $S_{2}^{\mathrm{rad}}(\lambda)=b_{\lambda}$. Since $b_{\lambda}=(\Lambda-\lambda)\hardy_{n}$ and $\Lambda=n^{2}/4$, (\ref{eq:S2}) follows. Moreover, since the Hardy inequality admits no extremal function, the same holds for $S_{2}(\lambda)$ and $S_{2}^{\mathrm{rad}}(\lambda)$.
\end{remark}

\section*{Acknowledgments}
The author would like to thank Roberta Musina for her useful comments after reading a preliminary version of this manuscript. 

\label{References}

\end{document}